\documentclass[12pt]{amsart}
\newcommand{\PSbox}[3]{\mbox{\rule{0in}{#3}\includegraphics{#1}\hspace{#2}}}
\title[Random Elements of Finite Distributive Lattices]{Generating Random Elements \\ of Finite Distributive Lattices}

\author{James Propp}
\address{Department of Mathematics, MIT}
\email{propp@math.mit.edu}

\thanks{During the conduct of the research that led to this article,
the author was supported by NSA grant MDA904-92-H-3060,
NSF grant DMS 9206374 and a grant from the MIT Class of 1922.}

\date{March 3, 1997; final revision January 13, 1998.}

\begin{document}

\begin{abstract}
This survey article describes a method for choosing uniformly at random 
from any finite set whose objects can be viewed as constituting a 
distributive lattice.  The method is based on ideas of the author 
and David Wilson for using ``coupling from the past'' to remove 
initialization bias from Monte Carlo randomization.  The article 
describes several applications to specific kinds of combinatorial 
objects such as tilings, constrained lattice paths, and alternating-sign 
matrices.

\end{abstract}

\maketitle

\begin{center}
{\it This article is dedicated to Herbert Wilf \\
in honor of his sixty-fifth birthday.}
\end{center}

\section{Introduction}

Herb Wilf, in addition to having done
important work
on problems related to counting combinatorial objects,
has also done pioneering research
on algorithms for generating combinatorial objects
``at random''
(that is, generating an element of a finite combinatorial set
so that each element has the same probability of being generated
as any other);
see \cite{CW}, \cite{DW}, \cite{GNW1}, \cite{GNW2}, \cite{NW1}, \cite{NW2}, 
\cite{NW3}, \cite{NW4}, \cite{Wilf1}, \cite{Wilf2} and \cite{Wilf3}
for fruits of this research.

This survey article describes a recent advance in the area
of random generation,
with applications to plane partitions, domino tilings,
alternating sign matrices,
and many other sorts of combinatorial objects.
The algorithm is of the ``random walk'' or ``Monte Carlo'' variety,
but unlike many such algorithms
it does not have any initialization bias.
The heart of the algorithm is the method of coupling from the past
explored by David Wilson and myself
in a joint article \cite{PW}.
For the sake of readability and motivation,
I will start by focusing on the application of our method
to plane partitions.

A beautiful formula of MacMahon \cite{MacMahon}
says that the number of plane partitions of $n$
with at most $a$ rows, at most $b$ columns,
and no part exceeding $c$
(hereafter to be called ``$(a,b,c)$-partitions'')
is given by
$$\prod_{i=0}^{a-1} \prod_{j=0}^{b-1} \prod_{k=0}^{c-1}
\frac{i+j+k+2}{i+j+k+1}$$
(see \cite{Andrews} and \cite{Stanley} for definitions
of ordinary partitions and plane partitions,
and section 2 of \cite{CLP} for a fairly simple self-contained proof
of MacMahon's formula).

Note that in the case $c=1$,
a plane partition with no part exceeding 1
can be read as
the Ferrers diagram of an ordinary partition,
and MacMahon's formula devolves into the assertion
that the number of ordinary partitions of $n$
with at most $a$ parts and no part exceeding $b$
is given by the binomial coefficient $\frac{(a+b)!}{a!b!}$.
Indeed, it is easy to see that such partitions
correspond to lattice paths of length $a+b$
joining $(0,a)$ to $(b,0)$,
or equivalently, combinations of $a+b$ elements taken $a$ at a time.
In view of these correspondences,
it is easy to generate a random $(a,b,1)$-partition.

Just as a lattice path in the $a$ by $b$ rectangle
is a 1-complex
(made up of edges in a 2-dimensional grid)
with prescribed boundary (namely the pair of points $(0,a)$ and $(b,0)$),
an $(a,b,c)$-partition
corresponds to a 2-complex 
(made up of square 2-cells in a 3-dimensional grid)
whose boundary is a particular non-planar hexagon
(namely the one that goes from
$(a,0,0)$ to $(a,b,0)$ to $(0,b,0)$ to $(0,b,c)$ to
$(0,0,c)$ to $(a,0,c)$ to $(a,0,0)$ in cyclic order).
In the former case, one requires that 
the lattice path should have (minimal) length $a+b$;
in the latter case, one requires that
the surface spanning the hexagon should have (minimal) area $ab+bc+ac$.

In this paper I will describe an algorithm
for generating a random $(a,b,c)$-partition
with $a,b,c$ arbitrary.
This algorithm was used to generate Figure 1,
which shows a random $(32,32,32)$-partition,
or rather the spanning surface that it determines,
viewed from a point on the ray $x=y=z>0$;
the three different orientations of grid-squares in 3-space
are seen in projection as three different orientations
of rhombuses in 2-space.
It should be stressed that the size of the plane partition
(that is, the sum of the parts)
was not specified in advance;
it is a random variable
with expected value $(32)^3/2$.

Note the non-homogeneity of the picture:
there is an approximately circular central region
in which rhombuses of different orientations are mixed together,
surrounded by six regions in which
rhombuses of a single orientation predominate.
As is shown in \cite{CLP},
in a certain strong probabilistic sense 
the boundary of the central region
does indeed converge to a perfect circle when $n$ is large.
This fact was first conjectured
on the basis of pictures like Figure 1,
and all known proofs depend on steps whose clearest motivation
comes from ``knowing the answer in advance''.
Thus we see that there are phenomena
pertaining to random $(a,b,c)$-partitions
that would not have been easy to divine by pure theory,
and that an algorithm for generating such plane partitions randomly
can be a valuable tool
for discovering and investigating such phenomena experimentally.

\begin{figure}
\begin{center}
\PSbox{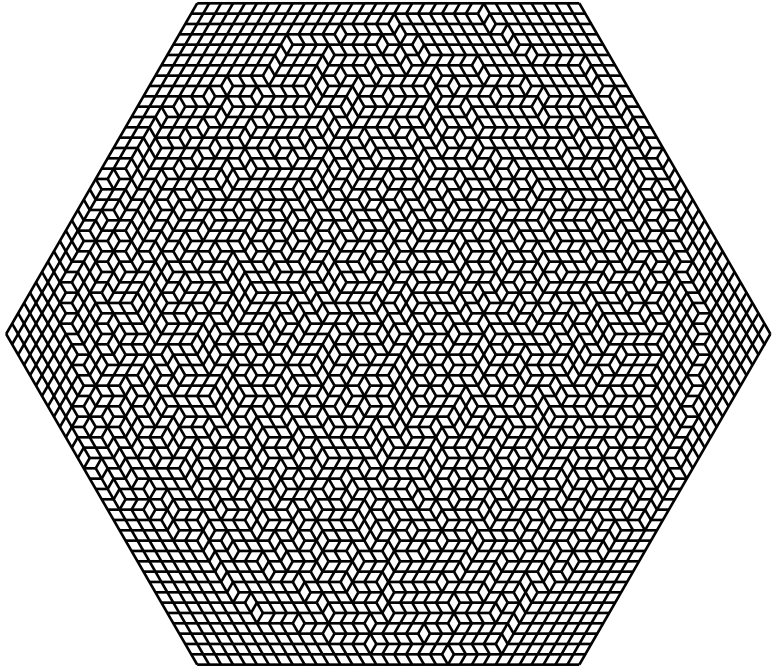}{288pt}{250pt}
\end{center}
\caption{A random lozenge tiling of a $32,32,32$ hexagon.}
\label{fig-random}
\end{figure}

In part 2 of this paper, I will show that this problem,
along with several others,
is a special case
of the problem of choosing a random element
of a finite distributive lattice.
In part 3, I will describe an algorithm
that allows one to solve this problem.
Like the well-known folk-algorithm
for generating a random element of a 3-element set
via independent tosses of an unbiased coin,
this algorithm runs in finite expected time,
even though it can take arbitrarily long to return an answer.
I conclude in part 4 with some open questions.

\section{Distributive lattices}

The 3-dimensional Ferrers diagram
% kludge
associated with an $(a,b,c)$-parti-
tion
is just an order ideal of the poset
obtained as the product of
chains of cardinalities $a$, $b$ and $c$.
Writing such chains as {\bf a}, {\bf b} and {\bf c}, respectively,
and using $J(\cdot)$ to represent the lattice of
order ideals of a finite poset (ordered by inclusion),
we see that what we are effectively trying to do
is generate a random element of
the distributive lattice $J({\bf a} \times {\bf b} \times {\bf c})$.
(See \cite{Stanley} for definitions
of posets, order ideals, lattices, and distributive lattices,
and for the fact that every finite distributive lattice
can be represented as $J(P)$ for some finite poset $P$.) 

Here is another way to view this ordering.
Consider a hexagon whose internal angles all measure 60 degrees
and whose sides in cyclic order have lengths $a,b,c,a,b,c$, respectively,
and dissect it into $2ab+2ac+2bc$ unit equilateral triangles.
Create a graph $G$ whose vertices correspond to these triangles,
where two vertices of $G$ are adjacent
if and only if
the corresponding triangles in the dissection share an edge.
($G$ is called a ``honeycomb graph'' for obvious reasons.)
Then each $(a,b,c)$-partition
corresponds to a tiling of the aforementioned hexagon
by rhombuses of side length 1
(each composed of two equilateral triangles joined edge to edge),
and each of these tilings in turn corresponds to 
a perfect matching of $G$,
that is, to a subset of the edges of $G$ in which
each vertex of $G$ appears just once.
It is shown in \cite{Propp}
(which was inspired by earlier work
of Conway and Lagarias \cite{CL}
and Thurston \cite{Thurston})
that the perfect matchings of any bipartite planar graph
can be given the structure of a distributive lattice.
In the case of the honeycomb-graph $G$ described above,
one obtains the same distributive lattice structure 
$J({\bf a} \times {\bf b} \times {\bf c})$ as before.  
If instead one allows graphs $G$ that are subgraphs of a square grid,
one is able to put a distributive lattice structure
on the set of all tilings of a 
connected subset of the square grid by dominos.
(These special cases are discussed in \cite{Thurston}.)

More generally (again see \cite{Propp}),
if $G$ is any bipartite planar graph,
and $d(\cdot)$ is any function from the vertex set of $G$
to the non-negative integers,
then the set of $d$-factors of $G$
(defined as the set of subgraphs of $G$
in which vertex $v$ has degree $d(v)$ for all vertices $v$)
can be given the structure of a distributive lattice
(outside of the trivial case in which no $d$-factors exist).

Another class of examples of distributive lattices in combinatorics
comes from constrained lattice paths.
For instance,
lattice paths of length $2n$
that go from $(0,n)$ to $(n,0)$
without straying outside of the triangle with vertices
$(0,0)$, $(n,0)$ and $(n,0)$
are a standard incarnation of the ``Catalan objects''.
More generally, we can consider all minimal-length lattice paths 
that go from one fixed lattice point to another 
without straying outside some fixed convex region.
Assuming that the set of such lattice paths is non-empty,
it is easy to show that the set of all such lattice paths
is a distributive lattice,
with the ordering being defined by inclusion of
the associated Ferrers diagrams.
One can also extend this approach
to the lattice paths enumerated by trinomial coefficients
(with step-vectors $(1,1)$, $(1,0)$, and $(1,-1)$).

A further example is provided by the alternating sign matrices
introduced by Mills, Robbins, and Rumsey \cite{MRR}.
One can transform each $n$ by $n$ alternating sign matrix
into an $n+1$ by $n+1$ array of numbers
in such a fashion that if one takes
the component-wise minimum (or component-wise maximum)
of any two such arrays,
one gets another array 
associated with an $n$ by $n$ alternating sign matrix.
In this fashion, the set of alternating sign matrices of fixed size
becomes a distributive lattice.
(See \cite{EKLP} for a more detailed discussion of this.)

A final situation I will mention in which
a non-obvious distributive lattice structure exists
is the set of independent sets
in a bipartite graph $G$.
(This lattice structure has been noticed
in the literature
several times independently.)
%This fact seems to have first been noticed by 
%{\sc someone or other} \cite{Who}.
Here is the way to see it:
If we color the vertices of $G$ black and white
and write every independent set $S$ as $S_B \cup S_W$
where the vertices in $S_B$ and $S_W$ are
black and white respectively,
then we may define the meet and join of $S_B \cup S_W$ and $T_B \cup T_W$
as  $(S_B \cup T_B) \cup (S_W \cap T_W)$
and $(S_B \cap T_B) \cup (S_W \cup T_W)$,
respectively.

We recall here that every finite distributive lattice
admits a unique representation (up to isomorphism) as $J(P)$,
for $P$ some finite poset
(which in fact is easily constructible from the lattice
as the sub-poset of join-irreducibles).
In each of the cases described above,
the bijection between the combinatorial objects described
and order ideals of an associated poset $P$
is easily implemented on a computer.
Thus, if we can solve the problem
of generating a random element of $J(P)$ for
an arbitrary finite poset $P$,
we will have solved 
the problem of generating random $d$-factors of bipartite planar graphs,
the problem of generating random independent sets 
in a general bipartite graph,
and many other problems as well.

\section{The algorithm}

In this section I describe a Monte Carlo approach
to sampling from the uniform distribution
on the set $J(P)$,
using a Markov chain whose state are the elements of $J(P)$
and whose steady state distribution is
the uniform distribution on $J(P)$.
If one merely simulated this Markov chain
for a large but finite number of steps
in the ordinary way,
one would have a distribution that was
close to the steady-state (uniform) distribution,
but there would be some residual error
(``initialization bias'').
Later in this section I will explain how one can get rid of this bias
by effectively simulating the Markov chain
for {\it infinitely} many steps,
from time minus infinity to time zero
(the method of coupling from the past).

It is natural to make $J(P)$ into a graph $H$
by declaring two order ideals to be adjacent
iff their symmetric difference consists of
exactly one element of $P$.
It is easy to show that this graph is connected.
Moreover, we can associate with each element $x$ of $P$
a randomization move
that preserves the uniform distribution on $J(P)$
(the probability distribution that
assigns each order ideal probability $1/|J(P)|$):
Given an order ideal $I$,
toss a fair coin,
and if the coin comes up heads (resp.\ tails),
replace $I$ by $I \cup \{x\}$
(resp., $I \setminus \{x\}$),
unless the resulting subset of $P$
is not an order ideal,
in which case leave $I$ alone.

If one performs an infinite sequence of such randomization moves
in which each new {\it randomization site} $x$ is chosen independently
from the uniform distribution on $P$
(or more generally from any probability distribution 
that assigns positive probability to each element of $P$),
then this process is simply a stationary Markov chain
whose state space is the set of order ideals of $P$;
standard ideas from the theory of Markov chains
guarantee that the probability of any particular order ideal
being the current order ideal
converges to $1/|J(P)|$ as time goes to infinity.

This Markov chain
gives us a way, in principle,
of generating a random element of $J(P)$
that is as close to unbiased as we like
(i.e., whose governing distribution is
as close to uniform as we like).
However, in the absence of estimates
of the mixing time of the Markov chain,
it is not clear for how long a time the chain must be run
in order to drive the bias below
some predetermined amount deemed acceptable.
Moreover, we seek a way of generating samples that
has {\it no} bias
when truly random bits are used.

Fortunately,
there is a way around this problem.
The Markov chain that we have described
is monotone
in the sense that if we were to run
{\it two} instantiations of it in parallel,
using the same randomization sites and the same coin-tosses in both runs
but starting from different initial order ideals $I_1$ and $I_2$,
then, provided $I_1 \supseteq I_2$,
the order ideals $I'_1$ and $I'_2$ that result after $n$ steps
of joint randomization
must satisfy the relation $I'_1 \supseteq I'_2$.
In particular, if we had chosen $I_1$ and $I_2$
to be the empty order ideal $\hat{0}=\emptyset$
and the full order ideal $\hat{1}=P$, respectively,
and if we find after $n$ steps that $I'_1$ and $I'_2$ are {\it equal}
(call their common value $I'$),
then {\it every} run of the Markov chain for $n$ steps
using those randomization sites and coin-tosses
will put us in state $I'$,
regardless of the initial state $I$;
all such histories are ``squeezed'' between
the $\hat{0}$-history and the $\hat{1}$-history,
and since these coalesce over the course of the simulation,
so must all histories.

It turns out that
for most of the examples that arise in practice,
this kind of coalescence occurs fairly quickly.
This gives us
a way of doing a kind of ``backwards simulation''
of the Markov chain
that effectively lets us run the chain
for a huge number of steps
without (usually) having to perform
anywhere near the number of steps
required in a straightforward simulation.
Say, for instance, that 
we want to run the Markov chain 
on the state-space $J(P)$
for one million steps,
starting from initial state $I^*$,
using random updates from time $-1,000,000$ to time 0
(we will see shortly why it is convenient
to index time in this way).
Suppose that our Markov chain is sufficiently rapidly mixing
that over the course of a thousand steps,
it's fairly likely for the initial states $\hat{0}$ and $\hat{1}$
to lead to the same final state
when evolved in tandem. 
If we simulate the Markov chain from time
$-1000$ to time 0,
using both $\hat{0}$ and $\hat{1}$ as initial states 
(let us call this ``phase one'' of our backwards simulation),
and we indeed find that both histories
coalesce at some state $I'$ at time 0,
then we do not need to run the full simulation
from one million steps in the past,
for we can already be sure that
such a simulation would have given us $I'$ as our sample.
Indeed, in this case we do not even have to choose
what the randomization sites and coin-tosses before time
$-1000$ are,
since they do not enter into the simulation.
In the unlikely event
that the two histories that were started at time $-1000$
do not coalesce by time 0,
then we could go back and do an honest simulation
for a million steps, starting from $I^*$ at time $-1,000,000$
(let us call this ``phase two'').
In this way we can simulate the behavior of
a million-step random walk
using (most of the time) only about two thousand simulation steps.
Note however that if one wishes to avoid
introducing bias into one's sample,
it is crucial that one use the {\it same} randomization sites
and coin-tosses from time $-1000$ onward
during the long, ``honest'' run (phase two)
as one did during the short preliminary run (phase one).

If one wants to simulate a {\it billion} steps of random walk,
one can modify phase two of the preceding algorithm by
running the million-step simulation 
(in the rare case where one thousand steps do not suffice)
using $\hat{0}$ and $\hat{1}$ as starting states,
rather than $I^*$.
Only in the incredibly rare case where this million-step simulation
fails to coalesce by time zero
would one need to resort to a billion-step simulation
starting from $I^*$ (``phase three''),
being careful once more to use the already-determined
randomization sites and coin-tosses from time $-1,000,000$ to time 0.

Note that the numbers one million and one billion
enter into these procedures in the form of bounds
on how far back into the past one is willing to go
before one ``honestly'' uses $I^*$ as the starting state
rather than trying to be clever by starting from $\hat{0}$ and $\hat{1}$
and hoping that they will coalesce by time 0.
Suppose that one removes this upper bound
on how far into the past one is willing to go:
if coalescence fails,
one goes back into the past 1000 times as far as one just did
and tries again.
Then one can show that with probability 1
the desired coalescence {\it will} eventually occur
(where our notion of ``eventuality'' goes backward in time
rather than forward),
and that the sample returned is effectively
a sample ``generated by a run of length infinity'' ---
that is, an absolutely unbiased sample.
This is David Wilson's method of coupling from the past.

It cannot be overemphasized that
the idea of progressing backwards into the past
is an indispensable feature of the algorithm.
In particular, if one were simply to run the Markov chain forwards
from initial states $\hat{0}$ and $\hat{1}$ in tandem 
until the histories coalesced
and then to output as one's sample the coalescent state,
one would in general get a biased sample.
In contrast,
samples generated via coupling from the past
are entirely free of bias,
to the extent that the coin-tosses used are random.

The version of the algorithm discussed above can be improved upon;
for instance, it is much better
to go progressively 1, 2, 4, 8, ... steps into the past
rather than $10^3$, $10^6$, $10^9$, ... steps.
Moreover, the algorithm can be generalized so as to apply
to interesting situations in statistical mechanics
where the desired distribution
is not uniform,
but is the Boltzmann distribution
for some non-constant energy function
on the configuration space.
All of these developments
are described more fully in \cite{PW}.

Here, we content ourselves with showing
that the expected running time is finite,
in the case where there exists a finite $L$ and a positive $\epsilon$
such that over the course of any time interval of length $L$ in the simulation,
the conditional probability of coalescence occurring,
given initial states $\hat{0}$ and $\hat{1}$,
is at least $\epsilon$.
(This holds for all the interesting applications;
in each case, it suffices to find an $L$ such that
the probability of going from $\hat{0}$ to $\hat{1}$ in $L$ steps, 
conditional upon starting in state $\hat{0}$,
is positive,
though the $L$ one gets from this approach is much, much larger
than the typical time-scale over which coalescence occurs.) 
To prove the claim,
note that over the course of $kL$ consecutive steps,
divided into $k$ blocks of length $L$,
the only way in which
coalescence could fail to occur
is if all $k$ blocks are ``non-coalescing'';
yet independence of the coin-tosses tells us
that this is an event of probability at most
$(1-\epsilon)^n$.  
As $n$ gets large, this probability shrinks to zero exponentially,
implying our claim on the finiteness of expected running time.

The algorithm described above is not always efficient;
indeed, in \cite{PW} Wilson and I describe a case in which
the modified Monte Carlo algorithm takes exponentially long,
even though a direct algorithm
for constructing a random order ideal
is quite easy to fashion
(see the Cautionary Note at the end of section 3.2 in \cite{PW}).
However, as a practical matter
the modified Monte Carlo algorithm is quite efficient
for many sorts of combinatorial objects
that arise naturally,
such as those mentioned in Section 2.

Three final points on implementation deserve note.

First, the algorithm remains valid if
instead of choosing randomization sites
in accordance with some fixed distribution on $P$
one chooses randomization sites according to some other
scheme, provided that two conditions are satisfied:
the choice of randomization sites
should not be affected by the outcomes of the coin flips,
and it must be the case that with probability 1
each element of $P$ gets chosen as a randomization site
infinitely often.
Under these hypotheses,
the (non-stationary) Markov chain on $J(P)$
has the uniform distribution as its
unique steady-state measure,
and the method of coupling from the past will get you there
with a finite amount of simulation
(with probability 1).
For instance, one may rotate among all the elements of $P$
in some fixed order,
rather than choosing the randomization sites randomly.
In this way one reduces the amount of information
that the algorithm needs to save.

Second, if one is using pseudo-random bits
given by some trusted pseudo-random number generator
(as I imagine most users of this algorithm will do),
then, though it is often necessary to {\it reuse} bits,
it is not necessary to {\it save} them all.
One can for instance save only seeds
that will enable one to re-create the formerly used bits
when they are needed again.

Third, it should be borne in mind that
the output of the Monte Carlo method
is likely to be somewhat correlated
to its running-time.
If one were to use the method
to generate $N$ samples,
where $N$ was determined on-the-fly
rather than chosen in advance,
then one might be contaminating one's samples with bias.
The scrupulous investigator may therefore wish
to commit to a certain value of $N$ ahead of time,
based on a preliminary investigation
of how long it is likely to take
to run the procedure $N$ times.

\section{Open problems}

Another sort of scheme for randomly generating an order ideal
of a finite poset $P$ is the recursive approach in which one 
decides whether or not to include $x$ in the random order ideal
by tossing a coin whose bias corresponds to the ratio
of the cardinalities of two distributive lattices,
namely, the lattice of order ideals of $P$ containing $x$
and the lattice of order ideals of $P$ not containing $x$.
This method (like coupling from the past) will not always be fast,
since the problem of counting antichains
(and hence, equivalently, order ideals)
in a finite partially ordered set
is \#P-complete (see \cite{PB}).
Yet another approach would be one along the lines of
recent work of Flajolet, Zimmermann, and van Cutsem
(see \cite{FZC1} and \cite{FZC2}),
who have already created Maple software
that (suitably instructed) can create
random plane partitions of the sort shown in Figure 1.
{\it Which way is ``best''?\/}
(Or rather: are there features of a lattice
that might dictate when one or another of these approaches
will do best?)

In the case of planar graphs $G$ that are not bipartite,
it is still possible to generate random perfect matchings of $G$
efficiently:
The method of Kasteleyn \cite{Kasteleyn}
allows one to write the number of perfect matchings
as a Pfaffian.
By calculating the ratio of two such Pfaffians,
one can determine the exact proportion of matchings of $G$
that contain a given edge,
and one can accordingly make an unbiased decision
as to whether or not to include this edge in the matching.
Applying this recursively,
one can generate a random matching of $G$.
Surprisingly, Wilson \cite{Wilson} showed that it is possible to randomly
generate perfect matchings within a constant multiple of the
time needed to compute just one Pfaffian; using sparse linear
algebra, this time is $O(n^{3/2})$ arithmetic operations,
or $O(n^{5/2} \:\, \log^2 n \:\, \log \log n)$ bit operations.
In the case where $G$ is bipartite,
a specialized version of the Pfaffian method
(the permanent-determinant method) applies
in much the same way;
but in this situation one can also
apply the Markov chain algorithm described in this paper.
{\it Can the method of this paper be extended
to apply in the non-bipartite case?}

If one uses a biased coin in place of a fair one,
one can devise an algorithm in which
the probability of a particular order ideal $I$ being generated
is proportional to $q^{|I|}$,
where $|I|$ is the cardinality of $I$
and $q$ is any positive real number.
That is, we can generate a random element of a distributive lattice
so that all the elements of rank $k$ in the lattice
have individual probability proportional to $q^k$ of being picked. 
We might try to choose $q$ so as to maximize the collective probability
of the elements of rank $k$,
but even using this optimal $q$
we might have to run the biased procedure many times
before we obtain a sample in rank $k$.
Note, however, that the resulting sample
will be governed by a uniform distribution
on the $k$th rank of the lattice.
{\it Is there an efficient procedure
for selecting an element in a particular rank
of a finite distributive lattice?}
Note that work of Flajolet, Zimmermann and van Cutsem  (\cite{FZC2})
provides a solution in some cases.

Lastly:
{\it To what extent can the method of coupling from the past
be extended to sampling from more general
partially ordered sets,\/}
such as modular lattices
or lattices in general?
For an example of an application of the method
to a special kind of non-distributive lattice,
see the paragraphs in subsection 3.3 of \cite{PW}
that treat permutations.

\bigskip

I thank the referee for many helpful comments.

\bigskip

{\sc Comments added after publication:}

The article gives the impression that the method of Section 3 was
used to generate Figure 1.  In fact, a more efficient approach
was used, exploiting the special structure of the lattice $J(P)$ in
question.  In particular, the poset $P$ of join-irreducibles of the 
lattice is a cube, and this cube can be divided into ``filaments''
(where elements  $...,(i,j,k),(i+1,j+1,k+1),...$ belong to the 
same filament).  Rather than merely attempting to add or delete 
a particular element of $P$, the algorithm attempts to modify $I$
by either deleting the largest element of the filament that is
in $I$ or adjoining the smallest element of the filament that is
not in $I$.

Putting it differently: the basic steps of the Markov chain are 
``hexagon-moves'' on the tiling, where one executes a non-trivial 
hexagon-move by changing the way a small unit hexagon of side-length 
1 inside the large region is tiled by three rhombuses (there are 
exactly two ways such a hexagon can be tiled).

Section 3 specifies that one is to choose the successive randomization
sites $x$ from the uniform distribution on $P$ (or more generally from any 
probability distribution that assigns positive probability to each element 
of $P$).  However, one has a great deal of leeway here.  For instance, one
can use a deterministic sequence of $x$'s that cycles through the set of
elements of $P$; as long as the coins used in the algorithm are random,
the procedure will still lead to a random order ideal.  In the case
where $P$ is ranked, one particularly natural idea is to do all the
$x$'s of even rank, followed by all the $x$'s of odd rank, and so on,
in alternation.  This particular scheme is well-suited to parallel
computation, since all the randomizations of the same parity that
are considered in any cycle can be done independently of one another.


\bigskip

\begin{thebibliography}{GNW2}

\bibitem[A]{Andrews}
G.\ Andrews,
{\sl The Theory of Partitions,\/}
Addison-Wesley, 1976.

\bibitem[CW]{CW}
E.\ Calabi and H.\ Wilf,
The sequential and random selection of subspaces over a finite field, 
{\it J.\ Combinatorial Theory\/} {\bf 22} (1977), 107--109.

\bibitem[CL]{CL}
H.\ Conway and J.\ Lagarias,
Tilings with polyominoes and combinatorial group theory,
{\it J.\ Combin.\ Theory A\/} {\bf 53} (1990), 183--208.


\bibitem[CLP]{CLP}
H.\ Cohn, M.\ Larsen and J.\ Propp,
The Shape of a Typical Boxed Plane Partition,
preprint, 1996.

\bibitem[DW]{DW}
J.\ Dixon and H.\ Wilf,
The random selection of unlabeled graphs, 
{\it J. Algorithms\/} {\bf 4} (1983), 205--213.

\bibitem[EKLP]{EKLP} 
N.\ Elkies, G.\ Kuperberg, M.\ Larsen and J.\ Propp, 
Alternating Sign Matrices and Domino Tilings, 
{\it J.\ Alg\  Combinatorics\/} {\bf 1} (1992), 
111--132 and 219--234.

\bibitem[FZC1]{FZC1}
P.\ Flajolet, P.\ Zimmermann and B.\ van Cutsem,
A calculus for the random generation of labelled combinatorial structures,
{\it Theoret.\ Comput.\ Sci.\/} {\bf 132} (1994), 1--35.

\bibitem[FZC2]{FZC2}
P.\ Flajolet, P.\ Zimmermann and B.\ van Cutsem,
A calculus of random generation: unlabelled structures,
preprint, 1996.

\bibitem[GNW1]{GNW1}
C.\ Greene, H.\ Nijenhuis and H.\ Wilf,
A probabilistic proof of a formula for the number of Young 
Tableaux of a given shape, 
{\it Adv.\ in Math.\/} {\bf 31} (1979), 104--109.

\bibitem[GNW2]{GNW2}
C.\ Greene, A.\ Nijenhuis and H.\ Wilf,
Another probabilistic method in the theory of Young Tableaux, 
{\it J.\ Combinatorial Theory\/} {\bf 37} (1984), 127--135.

\bibitem[K]{Kasteleyn}
P.W.\ Kasteleyn,
The statistics of dimers on a lattice, I. The number of dimer
arrangements on a quadratic lattice,
{\it Physica\/} {\bf 27},
1209--1225 (1961).

\bibitem[M]{MacMahon}
P.A.\ MacMahon, {\sl Combinatory Analysis\/},
Cambridge University Press, 1915--16 (reprinted by Chelsea Publishing
Company, New York, 1960).

\bibitem[MRR]{MRR}
W.\ Mills, D.\ Robbins, and H.\ Rumsey, Jr.,
Alternating sign matrices and descending plane partitions,
{\it J.\ Comb.\ Theory A\/} {\bf 34} (1983), 340--359.

\bibitem[NW1]{NW1}
A.\ Nijenhuis and H.\ Wilf,
A method and two algorithms in the theory of partitions, 
{\it J.\ Combinatorial Theory\/} {\bf 18} (1975), 219--222.
       
\bibitem[NW2]{NW2}
A.\ Nijenhuis and H.\ Wilf,
{\sl Combinatorial Algorithms\/}, Academic Press, 1975.

\bibitem[NW3]{NW3}       
A.\ Nijenhuis and H.\ Wilf,
{\sl Combinatorial Algorithms for Computers and Calculators\/} 
(second edition of \cite{NW2}), 
Academic Press, New York, 1978.  
       
\bibitem[NW4]{NW4}
A.\ Nijenhuis and H.\ Wilf,
The enumeration of connected graphs and linked diagrams, 
{\it J.\ Combinatorial Theory\/} {\bf 27} (1979), 356--359.

\bibitem[P]{Propp}
J.\ Propp, 
Lattice structure for orientations of graphs,
preprint, 1993.

\bibitem[PB]{PB}
J.\ Provan and M.\ Ball,
The complexity of counting cuts and of computing the probability that
a graph is connected,
{\it SIAM J.\ Comput.\/} {\bf 12} (1983), 777--788.

\bibitem[PW]{PW}
J.\ Propp and D.\ Wilson,
Exact sampling with coupled Markov chains and applications
to statistical mechanics,
{\it Random Structure and Algorithms\/}, to appear.

\bibitem[S]{Stanley}
R.\ Stanley, {\sl Enumerative Combinatorics I},
Wadsworth, 1986.

\bibitem[T]{Thurston}
W.\ Thurston, Conway's tiling groups,
{\it American Mathematical Monthly\/} {\bf 97} (1990), 757--773.

\bibitem[W1]{Wilf1}
H.\ Wilf,
A unified setting for sequencing, ranking and random 
selection of combinatorial objects, 
{\it Adv.\ in Math.\/} {\bf 24} (1977), 281--291.  
       
\bibitem[W2]{Wilf2}       
H.\ Wilf,
A unified setting for selection algorithms, II, 
{\it Annals of Discrete Mathematics\/} {\bf  2}: 
{\it Algorithmic aspects of combinatorics}, North Holland (1978), 135--148.
% Fix this reference.
  
\bibitem[W3]{Wilf3}       
H.\ Wilf,
The uniform selection of free trees, 
{\it J.\ Algorithms\/} {\bf 2} (1981), 204--207.

\bibitem[Wi]{Wilson}
D.\ Wilson,
Determinant algorithms for random planar structures,
to appear in the published abstracts of
the 1997 ACM-SIAM Symposium on Discrete Algorithms.

%\bibitem[Who]{Who}
%Who?  (I'm still tracking this reference down.)
%% Fix!

\end{thebibliography}
\end{document}